\documentclass[a4]{amsart}
\usepackage{setspace}
\usepackage{a4}
\usepackage{amssymb,amsmath,amsthm,latexsym}

\usepackage[left=20mm,top=0.6in,bottom=12mm]{geometry}

\usepackage{amsfonts}
\usepackage{amsfonts}
\usepackage{graphicx}
\usepackage{textcomp}
\usepackage{cite}
\usepackage{enumerate}
\usepackage[mathscr]{euscript}
\usepackage{mathtools}
\newtheorem{theorem}{Theorem}[section]

\newtheorem{definition}[theorem]{Definition}
\newtheorem{example}[theorem]{Example}

\newtheorem{remark}[theorem]{Remark}
\title{This is the title}
\usepackage{amssymb}
\usepackage{amssymb}
\usepackage{amssymb}
\usepackage{amssymb}
\usepackage{amsmath}
\usepackage{tikz}
\usepackage{hyperref}
\usepackage{enumerate}
\usepackage{mathtools}
\usepackage{amsmath}
\usepackage{tikz}
\usepackage{amssymb}
\usepackage{amsmath}
\usepackage{tikz}
\usepackage{hyperref}
\begin{document}
\title
{{P-Operators on Hilbert Spaces}}
\author{Rashid A. and P. Sam Johnson}
\address{Department of Mathematical and Computational Sciences, National Institute of Technology Karnataka (NITK), Surathkal, Mangaluru 575 025, India}
\email{rashid441188@gmail.com (First author) and sam@nitk.edu.in (Second author)}

\maketitle

\begin{abstract} 
A real square matrix $A$ is called a P-matrix if all its principal minors are positive. Using the sign non-reversal property of matrices, the notion of P-matrix has been recently extended by Kannan and Sivakumar to infinite-dimensional Banach spaces relative to a given Schauder basis. Motivated by their work,  we discuss P-operators on separable real Hilbert spaces. We also investigate P-operators relative to various orthonormal bases. \vspace*{5mm}

\end{abstract}



\section{Introduction}
An $n \times n$ real matrix $A$ is said to be a P-matrix \cite{MR142565} if all its principal minors are positive. The study of P-matrices originated in the context of some of the notable classes of matrices, such as positive matrices, M-matrices, and totally positive matrices. But the first systematic study of P-matrices appeared in the work of Fiedler and Pt\'{a}kk \cite{MR142565}. Since then, researchers have been captivated by this class of matrices. They play an important role in a wide range of applications, including the linear complementarity problem, global univalence of maps, linear differential inclusion problems, interval matrices, and computational complexity \cite{MR1298430, MR111,R3, R1,MR3396730,R2}. The linear complementarity problem (LCP) is stated as follows : Given an $n \times n$ real matrix $A$ and a vector $q \in \mathbb{R}^n$, the LCP is written as $LCP(A,q)$, and it is defined as to find a vector $x \in \mathbb{R}^n$ such that $x \geq 0$, $Ax+q \geq 0$ and $x^T(Ax+q) = 0$, where the notation $x \geq 0$ denotes each coordinate of the vector $x$ is non-negative. It is shown in \cite{MR3396730} that given a real square matrix $A$, the linear complementarity problem $LCP(A,q)$ has a unique solution for each  vector $q \in \mathbb{R}^n$ if and only if $A$ is a P-matrix.

 We say that an $n \times n$ matrix $A$ reverses the sign of a vector $x \in \mathbb{R}^n$ if $x_i(Ax)_i \leq 0$ for all $i=1,2,3, \ldots, n$, where $x_i$ denotes the $i^\text{th}$ coordinate of the vector $x$. Fiedler and Pt\'{a}kk \cite{MR142565} have shown that $A$ is a P-matrix if and only if $A$ does not reverse the sign of any non-zero vector. Inspired by this characterization of P-matrices, Kannan and Sivakumar \cite{MR3485836} extended the notion of P-operator to infinite-dimensional Banach spaces having a Schauder basis.  In this paper, we discuss the notion of P-operator to separable real Hilbert spaces and some results in this setting.
 
 In what follows, we will use separable real Hilbert space $\mathcal H$ and the term operator on $\mathcal H$ to mean a linear operator from $\mathcal H$ into itself. We denote $\mathcal B(\mathcal H)$ for the space of all bounded linear operators on $\mathcal H$. For $A\in \mathcal B(\mathcal H)$, we denote the adjoint of $A$ by $A^*$.

\section{P-Operators on Hilbert Spaces}
Let us begin with the definition of P-operator in Banach spaces introduced by Kannan and Sivakumar \cite{MR3485836}.  A  sequence $\{z_n\}_n$ in a real Banach space $X$ is said to be a Schauder basis for $X$ if for each $x\in X$, there exists a unique sequence of scalars $\{\alpha_n(x)\}_n$  such that $x=\sum_{n} \alpha_n(x) z_n.$  In such a case, we denote $x_n=\alpha_n(x)$, for any natural number $n$.  Throughout the set of natural numbers is the index set, and we write simply $\{\alpha_n(x)\}_n$ instead of $\{\alpha_n(x)\}_{n=1}^\infty.$

\begin{definition}\label{d2}\cite{MR3485836} 
Let $X$ be a Banach space with a Schauder basis. A bounded linear operator $T : X \rightarrow X$ is said to be a P-operator relative to the given Schauder
basis if for $x \in X$, the inequalities $x_n(Tx)_n \leq 0$ for all $n$ imply that $x = 0$.
\end{definition}

It is well-known that a countable orthonormal basis $\mathcal{B}=\{e_n\}_n$ exists for every separable Hilbert space $\mathcal{H}$ such that for any $x\in \mathcal{H}$ we have $x =\sum_n \langle x, e_n \rangle e_n$.  If an orthonormal basis is known, say $\{e_n\}_n$, then any orthonormal basis of $\mathcal H$ is of the form $\{Ue_n\}_n$ for some unitary operator $U$ on $\mathcal H$. We define the P-operator on separable real Hilbert spaces as follows.

\begin{definition}
Let  $\mathcal{B}=\{e_n\}_n$ be an orthonormal basis of $\mathcal{H}$. A bounded linear operator $T$ on $\mathcal{H}$ is said to be a P-operator relative to the given orthonormal basis $\mathcal{B}$ if for $x \in \mathcal{H}$, the inequalities $$\langle x, e_n \rangle \langle Tx, e_n \rangle \leq 0$$ for all $n$ imply that $x=0$.  
\end{definition}

\begin{example}\label{example-1}
Let $\ell_2$ denote the square summable sequence space of real numbers.  Let $\mathcal{B}= \{e_n\}_n$ be the standard orthonormal basis  of $\ell_2$, where $e_n$ denotes the vector whose $n^\text{th}$ entry is one, and all other entries are zero. Define $T: \ell_2 \rightarrow \ell_2$ by $$T(x_1,x_2,x_3,\ldots )= (\alpha_1x_1, \alpha_2x_2, \alpha_3x_3,\ldots),$$ for any $(x_1,x_2,x_3,\ldots)\in  \ell_2$ with $\alpha_n > 0$ for all $n$ and $\sup_n|\alpha_n|<\infty$. Then $T$ is a bounded linear operator, and it is a P-operator relative to $\mathcal{B}$. 
\end{example}

\begin{example}
The right shift operator $T_R$ and the left shift operator $T_L$ on  $\ell_2$ are defined by $$T_R(x_1, x_2, x_3, \ldots)=(0,x_1,x_2,\ldots)$$ and $$T_L(x_1, x_2, x_3, \ldots)=(x_2,x_3,x_4,\ldots)$$ respectively. The operators $T_R$ and $T_L$ are not P-operators relative to the standard orthonormal basis $\mathcal{B}= \{e_n\}_n$ of $\ell_2$.  Indeed, the non-zero element  $x=(1, -\frac{1}{2},\frac{1}{3}, -\frac{1}{4}, \ldots) \in \ell_2 $ satisfies the inequalities $\langle x, e_n \rangle \langle T_R(x), e_n \rangle \leq 0$  and $\langle x, e_n \rangle \langle T_L(x), e_n \rangle \leq 0$, for all $n$.
\end{example}

\begin{example}
The operators $I+T_R$ and $I+T_L$ on $\ell_2$ are
P-operators relative to the standard orthonormal basis $\mathcal{B}= \{e_n\}_n$ of $\ell_2$, where $I$ is the identity operator on $\ell_2$. Note that $I+T_R$ is a bounded linear operator. Suppose for $x=(x_1,x_2,x_3,\ldots ) \in  \ell_2$, the inequalities $\langle x, e_n \rangle \langle (I+T_R)x, e_n\rangle \leq 0$ for all $n$.  This leads to the inequalities $x_1^2 \leq 0$, $x_{n-1}x_n+x_n^2 \leq 0$, for all $n\geq 2$. From these inequalities, we get that  $x_n=0$, for all $n$, hence $x=0$.

Next, to see that $I+T_L$ is a P-operator relative to $\mathcal{B}$, it is noted that $I+T_L$ is a bounded linear operator. Suppose for $x =(x_1,x_2,x_3, \ldots ) \neq 0 \in \ell_2$ the inequalities $$\langle x, e_n \rangle \langle (I+T_L)x, e_n\rangle \leq 0$$ hold for all $n$.   This implies that  $x_n^2+x_nx_{n+1} \leq 0$ for all $n$.   Suppose $x_i=0$ for some index $i$, then we get that $x_{i-1}^2+x_{i-1}x_{i} \leq 0$, hence  $x_n=0$ for all $n \leq i$. Thus $x$ is of the form $x=(0,0, \ldots , 0, 0, x_{i+1}, \ldots )$ with $x_n \neq 0$ for all $n \geq i+1$ satisfying $x_n^2+x_nx_{n+1} \leq 0$.  Thus $|x_n^2| \leq |x_n|\ |x_{n+1}|$ for all $n\geq i+1$.

Now as $x_n \neq 0$ for all $n \geq i+1$, we have $|x_n|>0$ for all $n \geq i+1$. Thus by dividing the inequalities by $|x_n|$, we get $|x_n| \leq |x_{n+1}|$ for all $n \geq i+1$. This shows that the absolute values of the components of $x$ are increasing, and hence $n^\text{th}$ term of $x$ will not converge to $0$, hence $x \notin \ell_2$. Thus if $\langle x, e_n \rangle \langle (I+T_L)x, e_n\rangle \leq 0$ for all $n$, then $x$ must be equal to $0$.

\end{example}

We next give a result that says that an operator can be a P-operator relative to several orthonormal bases. 

\begin{theorem}\label{theorem-2}
Let $T$ be a bounded linear operator on $\mathcal{H}$ satisfying $TU=UT$  for a unitary operator $U$ on $\mathcal{H}$. Then $T$ is a P-operator relative to an orthonormal basis $\mathcal{B}=\{e_n\}_n$ of $\mathcal{H}$ if and only if $T$ is a P-operator relative to the orthonormal basis $\mathcal{B'}=\{Ue_n\}_n$ of $\mathcal{H}$.
\end{theorem}

\noindent {\bf Proof.} Suppose $T$ is a $P$-operator relative to the orthonormal basis $\mathcal{B}=\{e_n\}_n$ of $\mathcal{H}$ satisfying $TU=UT$.  Suppose $\langle x, Ue_n \rangle \langle Tx, Ue_n \rangle \leq 0$ for all $n$.  Then
$\langle U^*x, e_n \rangle \langle U^*Tx, e_n \rangle \leq 0$ for all $n$, hence $\langle U^*x, e_n \rangle \langle TU^*x, e_n \rangle \leq 0$ for all $n$, because  $TU=UT$.  As $T$ is a P-operator relative to the orthonormal basis $\mathcal{B}$, we get $U^*x=0$, hence $x=0$. Therefore $T$ is a P-operator relative to the orthonormal basis $\mathcal{B'}=\{Ue_n\}_n$.

\indent On the other hand, assume that $T$ is a P-operator relative to the orthonormal basis $\mathcal{B'}=\{Ue_n\}_n$ of $\mathcal{H}$ satisfying $TU=UT$. Suppose $\langle x, e_n \rangle \langle Tx, e_n \rangle \leq 0$ for all $n$.  As $U$ is a unitary operator, we get $\langle x, U^*Ue_n \rangle \langle Tx, U^*Ue_n \rangle \leq 0$ for all $n$, so $\langle Ux, Ue_n \rangle \langle UTx, Ue_n \rangle \leq 0$ for all $n$.  Hence $\langle Ux, Ue_n \rangle \langle TUx, Ue_n \rangle \leq 0$ for all $n$. Since $T$ is a P-operator relative to $\mathcal{B'}=\{Ue_n\}_n$, we get $Ux=0$, hence $x=0$. Therefore $T$ is a P-operator relative to the orthonormal basis $\mathcal{B}=\{e_n\}_n$.\hfill $\Box$\\

The condition $TU=UT$ in Theorem \ref{theorem-2} cannot be dropped which is illustrated in the example given below. The example also tells that an operator $T$ can be a P-operator relative to one orthonormal basis, whereas the same operator relative to another orthonormal basis may not be a P-operator. 

\begin{example}
Define $T: \ell_2 \rightarrow \ell_2$ by $$T(x_1, x_2,x_3, \ldots )= (x_1, 2  x_1+x_2, 2  x_2+x_3, \ldots )$$ for $x=(x_1,x_2, x_3, \ldots) \in \ell_2$. Then $T$ is a bounded linear operator and it is a P-operator relative to the standard orthonormal basis $\mathcal{B}=\{ e_n\}_n$ of $\ell_2$. 

Now consider the unitary operator $U$ on $\ell_2$ given by $$U(x_1,x_2,x_3, \ldots)=\Big(\frac{x_1}{\sqrt{2}}+\frac{x_2}{\sqrt{2}}, \frac{x_1}{\sqrt{2}}-\frac{x_2}{\sqrt{2}}, \frac{x_3}{\sqrt{2}}+\frac{x_4}{\sqrt{2}}, \frac{x_3}{\sqrt{2}}-\frac{x_4}{\sqrt{2}}, \ldots\Big),$$ then $U^* = U$ and $UT \neq TU$. The operator $T$ is not a P-operator relative to the orthonormal basis $\mathcal{B}^\prime=\{Ue_n\}$ of $\mathcal{H}$, because the non-zero element $x=(1, -1, 0, 0,\ldots) \in \ell_2$ satisfies the inequalities $\langle x, Ue_n \rangle \langle Tx, Ue_n \rangle \leq 0$,  for all $n$.
\end{example}

\begin{remark}
Theorem \ref{theorem-2} tells us that the condition $TU=UT$ is sufficient for the operator $T$ to be a P-operator relative to the orthonormal bases $\mathcal{B}=\{e_n\}_n$ and $\mathcal{B^\prime}=\{Ue_n\}_n$. But it is not a necessary condition. That is, an operator $T$ can be P-operator relative to two orthonormal bases $\mathcal{B}=\{e_n\}_n$ and $\mathcal{B^\prime}=\{Ue_n\}_n$, but it may not satisfy the relation $TU=UT$. The following example shows this fact.
\end{remark}

\begin{example}
Define $T: \ell_2 \rightarrow \ell_2$ by $$T(x_1, x_2,x_3, \ldots )= (x_1-x_2, x_1+x_2,   x_3-x_4, x_3+x_4, \ldots )$$ for $x=(x_1,x_2, x_3, \ldots) \in \ell_2$. Then $T$ is a P-operator relative to the standard orthonormal basis $\mathcal{B}=\{ e_n\}_n$ of $\mathcal{H}$. To see this, here the operator $T$ is bounded linear.  Suppose $\langle x, e_n \rangle \langle Tx, e_n \rangle \leq 0$ for all $n$.  Then $x_n(x_n-x_{n+1})\leq 0$  for odd $n$ and $x_{n}(x_{n}+x_{n-1})\leq 0$ for even $n$. Solving these inequalities together will lead to $x=0$.  Hence $T$ is a P-operator relative to $\mathcal{B}$.

Now consider the unitary operator $U$ on $\mathcal{H}$ given by $$U(x_1,x_2,x_3, \ldots)=\Big(\frac{x_1}{\sqrt{2}}+\frac{x_2}{\sqrt{2}}, \frac{x_1}{\sqrt{2}}-\frac{x_2}{\sqrt{2}}, \frac{x_3}{\sqrt{2}}+\frac{x_4}{\sqrt{2}}, \frac{x_3}{\sqrt{2}}-\frac{x_4}{\sqrt{2}}, \ldots\Big),$$ then $U^* = U$ and $\mathcal{B}^\prime=\{Ue_n\}_n$ is the another orthonormal basis of $\mathcal{H}$. Then $T$ is also a P-operator relative to $\mathcal{B}'$. To see this, suppose $\langle x, Ue_n \rangle \langle Tx, Ue_n \rangle \leq 0$ for all $n$. Then $x_n(x_n-x_{n-1})\leq 0$ for odd $n$ and $x_{n}(x_{n}+x_{n+1})\leq 0$ for even $n$.  Solving these inequalities together will lead to  $x=0$.  Hence $T$ is a P-operator relative to $\mathcal{B}'$. Note that  $$UT(x)=(\sqrt{2}x_1, -\sqrt{2}x_2, \sqrt{2}x_3, \ldots)$$ and $$TU(x)=(x_1,x_2,x_3,\ldots)$$ for any $x=(x_1,x_2,x_3,\ldots) \in \ell_2$, thus $TU \neq UT$.
\end{example}

\begin{theorem}

Let $\mathcal{B}=\{e_n\}_n$ be an orthonormal basis  of $\mathcal{H}$. Then the following statements hold good:\vspace{.3cm}
\indent (a) $T$ is a P-operator on $\mathcal{H}$ relative to $\mathcal{B}$ if and only if the operator $UTU^*$ is a P-operator relative to $\mathcal{B'}=\{Ue_n\}_n$, for any unitary operator U. \\
\indent (b)  $T$ is a P-operator on $\mathcal{H}$ relative to $\mathcal{B'}=\{Ue_n\}_n$ of $\mathcal{H}$ where $U$ is any unitary operator on $\mathcal{H}$ if and only if the operator $U^*TU$ is a P-operator on $\mathcal{H}$ relative to $\mathcal{B}$.
\end{theorem}

\noindent {\bf Proof.} (a) Assume that $T$ is a P-operator relative to $\mathcal{B}$. Then we have, if $\langle x, e_n \rangle \langle Tx, e_n\rangle \leq 0$ for all $n$ imply that $x=0$. Suppose $\langle x, Ue_n \rangle \langle UTU^*x, Ue_n\rangle \leq 0$ for all $n$.  Then $$\langle U^*x, e_n \rangle \langle TU^*x, e_n\rangle \leq 0$$ for all $n$.  As $T$ is a P-operator relative to $\mathcal{B}$, we get here $U^*x=0$ and hence $x=0$. Therefore $UTU^*$ is a P-operator relative to $\mathcal{B'}=\{Ue_n\}_n$.\\
\indent Conversely, assume that  $UTU^*$ is a  P-operator relative to the orthonormal basis $\mathcal{B'}=\{Ue_n\}_n$.  Suppose $\langle x, e_n \rangle \langle Tx, e_n\rangle \leq 0$ for all $n$.  Then  $\langle Ux, Ue_n \rangle \langle UTU^*Ux, Ue_n\rangle \leq 0$ for all $n$.  As $UTU^*$ is a P-operator relative to $\mathcal{B'}=\{Ue_n\}_n$, we get that $x=0$. Therefore $T$ is a P-operator relative to $\mathcal{B}$.   \\
\indent (b) Assume that $T$ is a P-operator relative to $\mathcal{B'}=\{Ue_n\}_n$. Suppose $\langle x, e_n \rangle \langle U^*TUx, e_n\rangle \leq 0$ for all $n$.  Then  $\langle Ux, Ue_n \rangle \langle TUx, Ue_n\rangle \leq 0$ for all $n$. As  $T$ is a P-operator relative to $\mathcal{B'}$, we get that $Ux=0$ and hence $x=0$. Therefore $U^*TU$ is a P-operator relative to $\mathcal{B}$.\\
\indent Conversely, assume that $U^*TU$ is a  P-operator relative to $\mathcal{B}$.  Suppose that  $$\langle x, Ue_n \rangle \langle Tx, Ue_n\rangle \leq 0$$ for all $n$.  Then $\langle U^*x, e_n \rangle \langle U^*TUU^*x, e_n\rangle \leq 0$ for all $n$. As $U^*TU$ is a P-operator relative to $\mathcal{B}$,  we get that $x=0$. Therefore $T$ is a P-operator relative to $\mathcal{B'}=\{Ue_n\}_n$. 
\hfill $\Box$\\

\begin{remark} 
	A bounded linear operator $T$ on $\mathcal H$ is	called invertible if there is a bounded linear operator $S$ on $\mathcal H$ so that $ST$ and $TS$ are the identity operators. We say that $S$ is the inverse of $T$ in this
	case and it is denoted by $T^{-1}$.
	It is observed in \cite{MR3485836} that every P-matrix is invertible and its inverse is also a P-matrix. But P-operator does not guarantee its invertibility in infinite-dimensional spaces which is shown in the following example.  Moreover, if we have a P-operator which is invertible, then its inverse is also a P-operator. 
\end{remark}

\begin{example}
	Consider the linear operator $T: \ell_2 \rightarrow \ell_2$ defined by $$T(x_1,x_2,x_3,\ldots)=\big(x_1,\frac{x_2}{2},\frac{x_3}{3},\ldots),$$ for $x=(x_1,x_2,x_3,\ldots\big)\in \ell_2$. Then T is a P-operator relative to the standard orthonormal basis of  $\ell_2$ but $T$ is not invertible.
\end{example}

\begin{theorem}\label{t14}
Let $T$ be an invertible P-operator on $\mathcal{H}$ relative to an orthonormal basis $\mathcal{B}= \{e_n\}_n$.  Then the inverse of $T$ is also a P-operator relative to $\mathcal{B}$.
\end{theorem}

\noindent {\bf Proof.}
Since $T$ is a P-operator on $\mathcal{H}$ relative to $\mathcal{B}$, we have, if the inequalities $\langle x, e_n \rangle \langle Tx, e_n \rangle \leq 0$ for all $n$ imply $x=0$. Suppose $\langle y, e_n \rangle \langle T^{-1}y, e_n \rangle \leq 0$ for all $n$. Hence $\langle Tx, e_n \rangle \langle x, e_n \rangle \leq 0$ for all $n$, where $x=T^{-1}y$. As $T$ is a P-operator,  we get $x=0$, hence $y=0$. Thus $T^{-1}$ is a P-operator relative to the orthonormal basis $\mathcal{B}$. 
\hfill $\Box$\\

\begin{theorem}\label{t15}
Let $T$ be a positive definite operator on  $\mathcal{H}$. Then $T$ is a P-operator on $\mathcal{H}$ relative to any orthonormal basis $\mathcal{B} = \{e_n\}_n$ of $\mathcal{H}$.  
\end{theorem}

\noindent {\bf Proof.} Assume that  $T$ is a positive definite operator on $\mathcal{H}$. Then for every $0 \neq x \in \mathcal{H}$, $\langle Tx, x \rangle > 0 $. Let $x$ be a non-zero element of $\mathcal{H}$. Then $ x=\sum_{n}x_ne_n $ and $Tx=\sum_n(Tx)_ne_n $. Then $\langle Tx, x \rangle = \sum_n(Tx)_nx_n  > 0$. Therefore there exists some $m$, for which $\langle x, e_m \rangle  \langle Tx, e_m \rangle > 0$. Hence if $\langle x, e_n \rangle  \langle Tx, e_n \rangle \leq 0$ for all $n$, imply that $x=0$. Hence $T$ is a P-operator. 
\hfill $\Box$\\

The converse of the above theorem need not be true which is shown in the following example.

\begin{example}\label{eg17}
 Let $T:\ell_2 \rightarrow \ell_2$ be defined by $$T(x_1, x_2, x_3, \ldots)= (x_1-7x_2, x_2, x_3, \ldots),$$ for $x=(x_1,x_2, x_3, \ldots) \in \ell_2$. 
  Then $T$ is a P-operator on $\ell_2$ relative to the standard basis $\mathcal{B}$ of $\ell_2$ as $\langle x, e_n\rangle \langle Tx, e_n \rangle \leq 0$ for all $n$ imply that $x=0$.  But $T$ is not a positive definite operator because $\langle Tx, x \rangle = x_1^2-7x_1x_2+x_2^2$ is negative for $x=(1,1, 0, \ldots)\neq 0$.
\end{example}

\subsection*{Acknowledgement}
The first author thanks the National Institute of Technology Karnataka (NITK) Surathkal for giving financial assistance.


\begin{thebibliography}{99}
\bibitem{MR1298430} Abraham Berman and Robert J. Plemmons, {\em Nonnegative matrices in the mathematical sciences},
 volume 9 of Classics in
Applied Mathematics, Society for Industrial and Applied Mathematics (SIAM), Philadelphia, PA, 1994.
\bibitem{MR1427262} Balmohan V. Limaye, {\em Functional analysis}, New Age International Publishers Limited,
New Delhi, second edition, 1996.


\bibitem{MR111}  Boyd S,  Ghaoui L. El,  Feron E and  Balakrishnan V, {\em Linear Matrix Inequalities in System
	and Control Theory}, SIAM,  Philadelphia, 1994.
	
\bibitem{R3} Christian Jansson and Jiri Rohn, {\em An algorithm for checking regularity of interval matrices}, SIAM
Journal on Matrix Analysis and Applications, 20, 756-776 (1999).

\bibitem{MR3121883} Melanie Lea Henthorn, {\em Examples of diagonal operators that fail spectral synthesis on spaces of analytic functions.} ProQuest LLC, Ann Arbor, MI, 2011, Thesis (Ph.D.)– Bowling Green State University.

\bibitem{R1} Miroslav Fiedler, {\em Special Matrices and Their Applications in Numerical Mathematics}, Martinus
Nijhoff Publishers, 1986.


\bibitem{MR142565} Miroslav Fiedler and Vlastimil Pt\'{a}kk, {\em On matrices with non-positive off-diagonal elements and positive principal minors},  Czechoslovak Math. J., {\bf 12(87)}, 382-400 (1962).



\bibitem{MR3485836}  Rajesh Kannan, M and  Sivakumar, K. C, {\em On certain positivity classes of operators}, Numer. Funct. Anal. Optim., {\bf37(2)}, 206--224 (2016).

\bibitem{MR3396730} Richard W. Cottle, Jong-Shi Pang, and Richard E. Stone, {\em The linear complementarity problem},  volume 60 of Classics
in Applied Mathematics, Society for Industrial and Applied Mathematics (SIAM), Philadelphia, PA, 2009.  

\bibitem{R2} Roger A. Horn and Charles R. Johnson, {\em Topics in Matrix Analysis}, Cambridge University Press, 1991.

\end{thebibliography}
\end{document}